\documentclass[11pt]{amsart2000}
\usepackage{amssymb}
\newtheorem{thm}{Theorem}

\newtheorem{prop}[thm]{Proposition}

\newtheorem{ex}[thm]{Example}
\theoremstyle{definition}
\newtheorem{defn}[thm]{Definition}

\newtheorem{qu}{Question}
\theoremstyle{remark}

\def\c{{\mathfrak c}}  \def\B{{\mathcal B}}
\def\C{{\mathcal C}}  \def\G{{\mathcal G}}
  
 \def\T{{\mathcal T}} 
 \def\W{{\mathcal W}}

\newcommand{\ran}{\operatorname{ran}}
\let\alp\alpha \let\bt\beta \let\gm\gamma \let\dlt\delta 
   \let\sgm\sigma 
 \let\w\omega \let\phi\varphi 
\let\nowt\varnothing
 
\newcommand{\dom}{\operatorname{dom}}

\begin{document}
\setlength{\unitlength}{0.01in}
\linethickness{0.01in}
\begin{center}
\begin{picture}(474,66)(0,0)
\multiput(0,66)(1,0){40}{\line(0,-1){24}}
\multiput(43,65)(1,-1){24}{\line(0,-1){40}}
\multiput(1,39)(1,-1){40}{\line(1,0){24}}
\multiput(70,2)(1,1){24}{\line(0,1){40}}
\multiput(72,0)(1,1){24}{\line(1,0){40}}
\multiput(97,66)(1,0){40}{\line(0,-1){40}}
\put(143,66){\makebox(0,0)[tl]{\footnotesize Proceedings of the Ninth Prague Topological Symposium}}
\put(143,50){\makebox(0,0)[tl]{\footnotesize Contributed papers from the symposium held in}}
\put(143,34){\makebox(0,0)[tl]{\footnotesize Prague, Czech Republic, August 19--25, 2001}}
\end{picture}
\end{center}
\vspace{0.25in}
\setcounter{page}{125}
\title{On the metrizability of spaces with a sharp base}
\author{Chris Good}
\address{School of Mathematics and Statistics, University of
Birmingham, Birmingham, B15 2TT, UK}
\email{c.good@bham.ac.uk}
\author{Robin W. Knight}
\address{Mathematical Institute, University of Oxford, 24-29 St Giles',
Oxford OX1 3LB, UK}
\email{knight@maths.ox.ac.uk}
\author{Abdul M. Mohamad}
\address{Department of Mathematics,
University of Auckland, Private Bag 92019, Auckland, New Zealand}
\email{mohamad@math.auckland.ac.nz}
\thanks{Reprinted from Topology and its Applications,
in press, 
Chris Good, Robin W. Knight and Abdul M. Mohamad,
On the metrizability of spaces with a sharp base,
Copyright (2002),
with permission from Elsevier Science \cite{gkm}.}
\thanks{Chris Good, Robin W. Knight and Abdul M. Mohamad,
{\em On the metrizability of spaces with a sharp base},
Proceedings of the Ninth Prague Topological Symposium, (Prague, 2001),
pp.~125--134, Topology Atlas, Toronto, 2002}
\begin{abstract}
A base $\mathcal{B}$ for a space $X$ is said to be {\it sharp} if, 
whenever $x\in X$ and $(B_n)_{n\in\w}$ is a sequence of pairwise
distinct elements of $\mathcal{B}$ each containing $x$,
the collection $\{\bigcap_{j\le n}B_j:n\in\w\}$ is a local base at $x$.
We answer questions raised by Alleche {\it et al.} and Arhangel'ski\u\i\
{\it et al.} by showing that a pseudocompact Tychonoff space with a sharp
base need not be metrizable and that the product of a space with a sharp
base and $[0,1]$ need not have a sharp base. We prove various metrization
theorems and provide a characterization along the lines of Ponomarev's
for point countable bases.
\end{abstract}
\maketitle

The notion of a uniform base was introduced by Alexandroff who proved that
a space (by which we mean $T_1$ topological space) is metrizable if and
only if it has a uniform base and is collectionwise normal \cite{al}. This
result follows from Bing's metrization theorem since a space has a uniform
base if and only if it is metacompact and developable. Recently Alleche,
Arhangel'ski\u\i\ and Calbrix \cite{AA} introduced the notions of sharp
base and weak development, which fit very naturally into the hierarchy of
such strong base conditions including weakly uniform bases (introduced by
Heath and Lindgren \cite{HL}) and point countable bases (see Figure 1
below). In this paper we look at the question of when a space, with a
sharp base is metrizable. In particular, we show that a pseudocompact
space with a sharp base need not be metrizable, but generalize various
situations where a space with a sharp base is seen to be metrizable.

\begin{defn} 
Let $\B$ be a base for a space $X$.
\begin{enumerate}
\item 
$\B$ is said to be {\it sharp} if, whenever $x\in X$ and $(B_n)_{n\in\w}$
is a sequence of pairwise distinct elements of $\B$ each containing $x$,
the collection $\{\bigcap_{j\le n}B_j:n\in\w\}$ is a local base at $x$.
\item 
$\B$ is said to be {\it uniform} if, whenever $x\in X$ and 
$(B_n)_{n\in\w}$ is a sequence of pairwise distinct elements of $\B$ each
containing $x$, then $(B_n)_{n\in\w}$ is a local base at $x$.
\item 
$\B$ is said to be {\it weakly uniform} if, whenever $\B^\prime$ is an 
infinite subset of $\B$, then $\bigcap\B^\prime$ contains at most one
point.
\item 
$\B$ is said to be a {\it weak development} if $\B=\bigcup_{n\in\w}\B_n$,
each $\B_n$ is a cover of $X$ and, whenever $x\in B_n\in\B_n$for each
$n\in\w$, then $\{\bigcap_{j\le n} B_j:n\in\w\}$ is a local base at $x$.
\end{enumerate}
\end{defn}

Arhangel'ski\u\i\ {\it et al.} prove that a space with a sharp base has a
point countable sharp base (\cite{AA} and \cite{A}) and is metaLindel\"of.
Moreover a weakly developable space has a $G_\dlt$-diagonal and a
submetacompact space with a base of countable order is developable
\cite{AA}.

We note in passing that the obvious definition of \lq uniform weak
developability\rq\ (having a base $\G=\bigcup\{\G_n:n\in\w\}$ such that
each $G_n$ is a cover and whenever $x\in G_n\in\G_n$, $\{G_n\}_n$ is a
base at $x$) is simply a restatement of developability. We also note that
a space with a $\sgm$-disjoint base need not have a sharp base: Bennett
and Lutzer \cite{BL2} construct a first countable (and a Lindel\"of)
example of a non-metrizable LOTS with $\sgm$-disjoint bases (and
continuous separating families), which can not have a sharp base by
Theorem \ref{metr}.

\begin{figure}[ht]
\begin{picture}(230,220)
\put(15,145){\makebox(0,0)[b]{\tiny development}}
\put(60,203){\makebox(0,0)[b]{\tiny uniform base $\equiv$}}
\put(60,192){\makebox(0,0)[b]{\tiny metacompact$+$development}}
\put(105,145){\makebox(0,0)[b]{\tiny sharp base}}
\put(137,147){\makebox(0,0)[b]{\tiny $\equiv$}}
\put(180,145){\makebox(0,0)[b]{\tiny point countable}}
\put(180,137){\makebox(0,0)[b]{\tiny sharp base}}
\put(60,40){\makebox(0,0)[b]{\tiny base of}}
\put(60,30){\makebox(0,0)[b]{\tiny countable order}}
\put(60,95){\makebox(0,0)[b]{\tiny weak}}
\put(60,85){\makebox(0,0)[b]{\tiny development}}
\put(150,95){\makebox(0,0)[b]{\tiny weakly}}
\put(150,87){\makebox(0,0)[b]{\tiny uniform base}}
\put(230,95){\makebox(0,0)[b]{\tiny point countable}}
\put(230,87){\makebox(0,0)[b]{\tiny base}}
\put(60,80){\vector(0,-1){30}}
\put(185,135){\vector(1,-1){30}}
\put(110,140){\vector(1,-1){35}}
\put(100,140){\vector(-1,-1){35}}
\put(65,190){\vector(1,-1){35}}
\put(55,190){\vector(-1,-1){35}}
\put(20,140){\vector(1,-1){35}}
\end{picture}
\caption{}
\end{figure}

When is a space with a sharp base metrizable?
We summarize the relevant results of \cite{AA}, \cite{A} and \cite{BL} in 
the following theorem.

\begin{thm}\label{metr}
Let $X$ be a regular space with a sharp base, then $X$ is metrizable if
any of the following hold:
\begin{enumerate}
\item $X$ is separable;
\item $X$ is locally compact (so a manifold with sharp base is metrizable);
\item $X$ is countably compact; \label{cc}
\item $X$ is pseudocompact and CCC;\label{ccc}
\item $X$ is a GO space.
\end{enumerate}
\end{thm}

A space is pseudocompact if every continuous real valued function is
bounded. Every (Tychonoff) pseudocompact space with a uniform base is
metrizable (see \cite{w}, \cite{sc} or \cite{us}), whilst a pseudocompact
space with a point-countable base need not be metrizable \cite{sha}.
Moreover pseudocompact Tychonoff spaces with regular $G_\dlt$-diagonals
are metrizable \cite{Mc}, whilst Mrowka's $\Psi$ space is an example of a
pseudocompact, non-metrizable Moore space. So it is natural to ask (see
\cite{AA} and \cite{A}) whether every pseudocompact space with a sharp
base is metrizable. The space $P$ of Example \ref{pseudex} shows that the
answer to this question is \lq no.\rq\ In addition, $P$ answers a number
of other questions in the negative: Alleche {\it et al.} ask whether the
product $X\times[0,1]$ has a sharp base if $X$ does; Heath and Lindgren
\cite{HL} ask whether a space with a weakly uniform base has a
$G^*_\dlt$-diagonal; and $P$ is another example (see \cite{sha} and
\cite{wat}) of a pseudocompact space with a point countable base that is
not compact, and is a non-compact pseudocompact space with a weakly
uniform base, answering questions of Peregudov \cite{per}.

\begin{ex}\label{pseudex}
There exists a Tychonoff, non-metrizable pseudocompact space with a sharp
base but without a $G_\dlt^*$-diagonal whose product with the closed unit
interval does not have a sharp base.
\end{ex}

\begin{proof} 
Our example is a modification of the example of a non-developable space
with a sharp base \cite {AA}. 
We add extra points to a (non-separable) metric space $B$ in such a way
that the resulting space is pseudocompact, has a sharp base but is not
compact, hence not metrizable.

Let $B={^\w\c}$ be the Tychonoff product of countably many copies of the
discrete space of size continuum with the usual Baire metric. 
For each finite partial function $f\in{^{<\w}\c}$, let $[f]$ denote the
basic open subset of $B$, $[f]=\{g\in{^\w\c}:f\subseteq g\}$ (so $[f]$ is
the collection of all elements of $B$ which agree with $f$ on 
$\dom f$). 
Note that, if $\dom f\subseteq\dom g$, then the two basic open sets $[f]$
and $[g]$ have non-empty intersection if and only if $f\subseteq g$ if and
only if $[g]\subseteq[f]$. 
If $[f]\cap[g]=\nowt$ then the functions $f$ and $g$ are incompatible (we
write $f\perp g$) and neither $f\subseteq g$ nor $g\subseteq f$.

Let 
$$\mathcal S=\{S\in{^\w({^{<\w}\c})}:\text{$S(m)\perp S(n)$,
for each $m$ and $n$}\},$$ 
so that each $S$ in $\mathcal S$ codes for a sequence of disjoint basic
open sets in $B$.
Enumerate $\mathcal S$ as $\{S_\alp:\alp\in\c\}$ in such a way that each
$S$ in $\mathcal S$ occurs $\c$ times. 
To ensure that our space is pseudocompact, we recursively add limit points
(to some of) these sequences of open sets. 
These limit points $s_\alp$ will have basic open neighbourhoods of the
form 
$$N(\alp,n)=\{s_\alp\}\cup\bigcup_{m\ge n}[T_\alp(m)],$$
where $T_\alp\in{^\w({^{<\w}\c})}$ is defined depending on $S_\alp$.

Suppose that for each $\alp<\gm$ we have either defined if possible a
sequence $T_\alp\in{^\w({^{<\w}\c})}$ such that
\begin{enumerate}
\item[($1\gm$)] 
for $i\neq j$, $T_\alp(i)\perp T_\alp(j)$,
\item[($2\gm$)] 
for $\bt<\gm$, $\bt\neq\alp$, $T_\bt$ defined, 
$\ran T_\alp\cap \ran T_\bt=\nowt$, and
\item[($3\gm$)] 
for $\bt<\gm$, $\bt\neq\alp$, $T_\bt$ defined, if 
$T_\alp(i)\supseteq T_\bt(j)$, then 
$T_\alp(i^\prime)\perp T_\bt(j^\prime)$ for all 
$\langle i^\prime,j^\prime\rangle\neq\langle i,j\rangle$
\end{enumerate}
or we have not defined $T_\alp$. We now define $T_\gm$.

First note that if $S_\gm^\prime(i)$ extends $S_\gm(i)$, then the open set
$[S_\gm^\prime(i)]$ is a subset of $[S_\gm(i)]$, so any limit of the
sequence of open sets $\{[S_\gm^\prime(i)]:i\in\w\}$ will also be a limit
of the sequence $\{[S_\gm(i)]:i\in\w\}$.

Since each $T_\alp(j)$ is finite, there is some $\dlt<\c$ which is not in
$\bigcup\{T_\alp(j):\alp<\gm, j\in\w\}$. 
For each $i\in\w$, let
$S^\prime_\gm(i)={S_\gm(i)}^\frown\{\dlt\}$ extend $S_\gm(i)$. Then for
all $i,j\in\w$ and $\alp<\gm$, $S^\prime_\gm(i)\nsubseteq T_\alp(j)$
and $T_\alp(j)\subseteq S^\prime(i)$ only if
$T_\alp(j)\subseteq S(i)$. Notice that this implies that
$[T_\alp(j)]\nsubseteq[S_\gm^\prime(i)]$ and that
$[S^\prime_\gm(i)]\subseteq[T_\alp(j)]$ only if
$[S_\gm(i)]\subseteq[T_\alp(j)]$.

Case 1: Suppose that there exists some $\alp<\gm$ for which $T_\alp$ was
defined,
such that for infinitely many $i\in\w$ there
exists some $j\in\w$ such that $S_\gm^\prime(i)\supseteq
S_\gm(i)\supseteq T_\alp(j)$. In this
case we do not define $T_\gm$ (since
infinitely many of the basic open sets $[T_\alp(j)]$ contain
an open set $[S_\gm(i)]$ and the limit point $s_\alp$ will deal
with the sequence $S_\gm$).

Case 2: Now suppose that Case 1 does not hold and that hence
\begin{enumerate}
\item[(*)]
for each $\alp<\gm$ there are at most finitely many
$i$ for which $S^\prime_\gm(i)\supseteq T_\alp(j)$ for some $j$.
\end{enumerate}
Suppose further that for each $i\le k$, we have chosen natural numbers
$0=r_0<r_1<\dots<r_k$ and defined
$T_\gm(i)$ to be $S^\prime_\gm(r_i)$.

Since each $T_\gm(i)$ is a finite partial function, there are at most finitely
many
possible partial functions such that $f\subseteq T_\gm(i)$ for some $i\le k$.
By condition ($2\gm$) there are at most finitely many
$\alp<\gm$ with such
an $f$ in $\ran T_\alp$. List these $\alp$ as $\alp(1),\dots\alp(m)$.
By (*), for each $\alp(m)$,
there is a $j_m$ such that for all $i\ge j$, $S_\gm^\prime(i)$ does not extend
any $T_{\alp(m)}(j)$. Now let $r_{k+1}=\max{j_m}$ and
$T_\gm(k+1)=S_\gm^\prime(r_{k+1})$.

We now claim that conditions ($1\c$), ($2\c$) and ($3\c$) hold.
Suppose that $T_\bt$ and $T_\alp$ were defined for some
$\bt<\alp<\c$. Condition ($1\c$) is obvious since each $T_\alp$ is
a subsequence of $S^\prime_\alp$ each term of which extends the
corresponding term of $S_\alp$, and $S_\alp$ is a sequence of
pairwise incompatible partial functions. ($2\c$) holds since, if
$\bt<\alp$, then the extension $S^\prime_\gm(i)$ was chosen to
ensure that $T_\bt(j)\nsupseteq S^\prime_\alp(i)$ for any $j$, so
in particular $T_\bt(j)\neq T_\alp(i)$ and $\ran T_\bt\cap\ran
T_\alp$. To see that ($3\c$) holds, note first that
$S_\alp^\prime(i)$ was chosen so that $S^\prime_\alp(i)\nsubseteq
T_\bt(j)$ for any $j$, which implies that $T_\alp(i)\nsubseteq
T_\bt(j)$ for any $\langle i,j\rangle$. On the other hand, suppose that 
$i$
is least such that for some $j$, $T_\bt(j)\subseteq T_\alp(i)$. If
$k>i$, then $T_\alp(k)=S_\alp^\prime(r_k)$ and $r_{k}$ was chosen
precisely so that $S_\alp^\prime(r_k)\nsupseteq T_\bt(l)$ for any
$l\in\w$. Moreover, there can be at most one $j$ such that
$T_\alp(i)\supseteq T_\bt(j)$, since by ($1\c$), $T_\bt(j)\perp
T_\bt(l)$, $j\neq l$. This completes the recursion.

Let $L=\{s_\alp:\text{$T_\alp$ has been defined}\}$ be a set of pairwise
distinct points disjoint form $B$ and let
$P=B\cup L$. We topologize $P$ by letting $B$ be an open subspace with the
usual
Baire metric topology and declaring the $n^{\text{th}}$
basic open set about the point $s_\alp$ to be the set
$N(\alp,n)=\{s_\alp\}\cup\bigcup_{m\ge n}[T_\alp(m)]$.

If $\T_\alp=\{[T_\alp(n)]:n\in\w\}$, then condition ($1\c$) ensures that each
$\T_\alp$
is a pairwise disjoint collection, ($2\c$)
ensures that each basic open set $[f]$ occurs in at most one $\T_\alp$,
and ($3\c$) ensures that if $N(\alp,n)$ meets $N(\bt,m)$, then
$N(\alp,n)\cap N(\bt,m)=[T_\alp(j)]\cap[T_\bt(k)]$ for some $j\ge n $ and
$k\ge m$.

That $P$ has a sharp base follows exactly as for the example due to
Alleche {\it et al}. 
Let $\B_B$ be a sharp base for $B$ and let
$$\B=\B_B\cup\{N(\alp,n):s_\alp\in L\text{ and }n\in\w\}.$$
Suppose $x\in\bigcap_{k\in\w} B_k$ for some (injective) sequence
$\{B_k\in\B:k\in\w\}$. Since $\B_B$ is a sharp base and $s_\alp\in N\in\B$ if
and
only if $N=(\alp,n)$ for some $n$, the only case that is not obvious is
when $x\in B$ and $B_k=N(\alp_k,m_k)$ for all but finitely many $k$.
But in this case condition ($3\c$)
implies that, for $n\ge1$, $\bigcap_{k\le n}B_k=\bigcap_{k\le
n}[T_{\alp_k}(j_k)]$. Moreover ($2\c$) implies that $T_{\alp_k}(j_k)\neq
T_{\alp_{k^\prime}}(j_{k^\prime})$, so that $\{\bigcap_{k\le
n}B_k:n\in\w\}$ contains a strictly decreasing subsequence and is
therefore a base at $x$.

Since the set $\{s_\alp:\alp\in\c\}$ is infinite, closed discrete, $P$ is not
compact.
On the other hand, $P$ is pseudocompact (so $P$ is not metrizable).
To see this, suppose that $\phi$ is a continuous real-valued function on $P$
taking
values in $[n,\infty)$
for each $n\in\w$. Since $B$ is
dense in $P$, for each $n\in\w$, there is some $x_n$ in $B$ such that
$\phi(x_n)>n$. By continuity, $\{x_n:n\in\w\}$ does not have a limit
point in $B$. Since $\phi$ is continuous and $B$ is metrizable, there are basic
open sets
$[f_n]$ for each $n\in\w$ such that $x_n\in[f_n]\subseteq \phi^{-1}(n,\infty)$
and $\{[f_n]:\in\w\}$
is a disjoint collection. But in this case $f_n\perp f_m$ when $n\neq m$
so that $\{f_n:n\in\w\}=S_\alp$ for some $\alp\in\c$. In which case,
either $s_\alp$ and $T_\alp$ were defined or $s_\alp$ was not defined
and, for some $\bt<\alp$, $T_\bt(j)\subseteq S_\alp(n)=f_n$ for
infinitely many $n$. In the second case, each basic open neighbourhood
$N(\bt,n)$
of $s_\bt$ contains infinitely many of the sets $[f_n]$. In the first case,
$T_\alp$
was chosen so that $T_\alp(i)\supseteq f_{r_i}$ for each $i\in\w$, so that
$[T_\alp(i)]\subseteq[f_{r_i}]$. In either case, each neighbourhood of
$s_\bt$ or $s_\alp$ contains points which take arbitrarily large values
under $\phi$, contradicting continuity.

Now suppose for a contradiction that $P\times [0,1]$ has a sharp base.
We shall show that this would imply that $P$ has a $\sgm$-point finite base,
which is impossible since Uspenski\u\i\ \cite{us} shows that a
pseudocompact space with a $\sgm$-point finite base is metrizable.

To this end, let $\W$ be a sharp base for $P\times [0,1]$
and let $\C$ be a countable sharp base for $[0,1]$.
For each $x$ in $L$ choose $W_n^x$ in $\W$,
$B_n^x$ in $\B$ (the sharp base for $P$), and $C_n^x$ in $\C$ such
that $B_n^x\times C_n^x\subseteq W_n^x$,
$\{W^x_n:n\in\w\}$ (and hence $\{B_n^x\times C_n^x:n\in\w\}$) is a local base at
$(x,1/2)$
and $W_0^x\cap(L\times[0,1])\subseteq \{x\}\times[0,1]$, which is possible since
$L$
is a closed discrete subset of $P$.

Let 
$$\B_C = \{B\in\B:\text{for some $n\in\w$ and some $x\in L$, $B=B_n^x$ and
$C=C^x_n$}\}.$$
If $\B_C$ is not point finite then
for some $y$ in $P$, $y\in\bigcap_{j\in\w}B_j$ for some pairwise distinct
$B_j\in\B_C$.
By definition, for each $j$ there is some $x_j\in L$ and
$n_j\in\w$ such that $B_j=B^{x_j}_{n_j}$ and $C=C^{x_j}_{n_j}$. But then
$$
\{y\}\times C \subseteq
\bigcap_{j\in\w}(B_{n_j}^{x_j}\times C_{n_j}^{x_j}) \subseteq
\bigcap_{j\in\w} W_{n_j}^{x_j}.
$$
Since $B_j\neq B_k$, either there is an infinite set $J\subseteq\w$ such that
$x_j\neq x_k$, for distinct $j,k\in J$, or there is an infinite set
$K\subseteq\w$ such that $x_j=x_k=x$ but $n_j\neq n_k$ for some $x\in L$ and
distinct $j,k\in K$. In the first case, $\{W_{n_j}^{x_j}:j\in J\}$ is a
pairwise
distinct subset of the sharp base $\W$ and $\bigcap_{j\in
J}W_{n_j}^{x_j}$ contains at most one point. In the second case
$\bigcap_{k\in K}(B_{n_k}^{x_k}\times C_{n_k}^{x_k})=(x,1/2)$, since
$\{B_n^x\times C_n^x:n\in\w\}$ is a local base at $(x,1/2)$. In either
case, $\{y\}\times C$ contains at most one point, which is not the case, and
$\B_C$ is point
finite.

Since $\{B_n^x\times C_n^x:n\in\w\}$ is a local base at $(x,1/2)$ and
$\C$ is countable, $\B=\bigcup_{C\in\C}\B_C$ is a $\sgm$-point finite
base for points of $L$. But $P=B\cup L]$
and $B$ is a metric space, so $P$ has a
$\sgm$-point finite base: a contradiction.

By Theorem \ref{psmoore}, $P$ does not have a $G_\dlt^*$ diagonal, nor
indeed is it submetacompact.
\end{proof}

So when is a pseudocompact space with a sharp base metrizable? 
As mentioned above, a pseudocompact, CCC regular space with a sharp base
is metrizable \cite[Theorem 21]{A}.
Pseudocompact, Moore spaces are CCC. 
Moreover, in proving that a pseudocompact Tychonoff space with a regular
$G_\dlt$-diagonal is metrizable, McArthur \cite{Mc} proves that a
pseudocompact space with a $G_\dlt^*$-diagonal is developable. 
Hence we have

\begin{thm}\label{psmoore}
A pseudocompact regular space $X$ with a sharp base is metrizable if
either of the following hold:
\begin{enumerate}
\item $X$ is developable, or;
\item $X$ has a $G^*_\dlt$-diagonal.
\end{enumerate}
\end{thm}

A pseudocompact space with a $G_\dlt$-diagonal is \v Cech complete
\cite[Lemma 20]{A}, hence Baire, so the following theorem is a
strengthening of Theorem 21 of \cite{A}.
A space is strongly quasi-complete if there is a map $g$ assigning to
each $x\in X$ and $n\in\w$ an open set $g(n,x)$ containing $x$ such that
$\{x_n\}$ clusters at $x$ whenever 
$\{x,x_n\}\subseteq \bigcap_{i\le n}g(i,y_i)$. Weakly developable spaces
are clearly strongly quasi-complete.

\begin{thm}\label{baire} 
A regular, locally CCC, locally Baire space with a sharp base is
metrizable.
\end{thm}

\begin{proof} 
Let $X$ be a regular, locally CCC, locally Baire
space with a sharp base. Since $X$ has a weak development, it is strongly
quasi-complete. Hodel \cite{H} shows that every regular, quasi-complete
CCC Baire space with either a $G_\dlt$-diagonal or a point countable
separating open cover is separable. Since $X$ has a sharp base, $X$ has
a point countable base, a $G_\dlt$-diagonal and is quasi-complete.
Hence $X$ is locally separable. But every locally separable regular space
with a point countable base is a disjoint union of clopen subspaces each of
which
has a countable base (see Theorem 7.2 of \cite{G}). Hence $X$ is metrizable.
\end{proof}

Generalising the fact that a countably compact space with a sharp base
is metrizable we have:

\begin{thm} A regular, $\w_1$-compact space with a sharp base is
metrizable. \end{thm}

\begin{proof} Since $X$ is $\w_1$-compact, every point-countable open cover of
$X$ has a countable subcover (Lemma 7.5, \cite{G}). Since $X$ has a sharp
base, it has a point countable base and therefore is Lindel\"of.
A metacompact space with a sharp base is developable \cite{AA} and so a
Lindel\"of space with a sharp base is metrizable.
\end{proof}

Not surprisingly a monotonically normal space with a sharp base is
metrizable (c.f.\ \cite{BL} where it is shown that a GO-space with
a sharp base is metrizable).

\begin{thm} \label{t} For a monotonically normal space $X$
the following are
equivalent:
\begin{enumerate}
\item $X$ is metrizable;
\item $X$ has a sharp base;
\item $X$ has a weak development;
\item $X$ is strongly quasi-complete;
\item $X$ has a base of countable order and a $G_\dlt$-diagonal.
\end{enumerate}
\end{thm}

\begin{proof} Since $1\implies2\implies3\implies4\implies5$
(that 4 implies 5 follows from Theorems 2.2 and 2.3 of \cite{gi}),
it remains to show
that a monotonically normal space with a base of countable order and a
$G_\dlt$-diagonal is metrizable. By the Balogh-Rudin theorem \cite{br}, since a
stationary
set of a regular cardinal does not have a $G_\dlt$-diagonal,
a monotonically
normal space with a $G_\dlt$-diagonal is paracompact. The result
then follows
since a paracompact space with a base of countable order is metrizable
\cite{arh}.
\end{proof}

The proof that $P\times[0,1]$ does not have a sharp base does not quite
extend to a proof that if the product of a space $X$ with $[0,1]$ has a
sharp base then $X$ has a $\sgm$-point finite base. The converse however
is easily seen to be true.

\begin{prop} 
If a space $X$ has a $\sgm$-point finite sharp base then $X\times[0,1]$
has a sharp base.\label{prod}
\end{prop}

\begin{proof} 
Suppose that $\B=\bigcup\B_n$ is a $\sgm$-point finite sharp base for $X$ 
and $\C=\bigcup\C_n$ is a development for $[0,1]$ such that each 
$\C_{n+1}$ is finite and refines $\C_n$ (so that $\C$ is also a sharp
base for $[0,1]$).
For each $n\in \w$ let $\W_n=\{B\times C:B\in\B_n, C\in\C_n\}$ and let
$\W=\bigcup_n\W_n$.

Firstly note that $\W$ is a base for $X\times[0,1]$. If $(x,r)$ is in
some open set $U$, choose $n$ and $B\in\B_m$ such that $(x,r)\in B\times
st(r,\C_n)\subseteq U$.
Now for some $k\ge\max\{m,n\}$, there is $B^\prime\in\B_k$ 
$x\in B^\prime\subseteq B$. 
But then, since $\C_k$ refines $\C_n$, if $r\in C\in\C_k$, 
$B^\prime\times C\in\W_k$ and
$$
(x,r)\in B^\prime \times C\subseteq
B^\prime\times st(r,\C_k)\subseteq 
B\times st(r,\C_n)\subset 
U.
$$

Now suppose that $(x,r)\in B_j\times C_j=W_j\in\W$ for distinct $W_j$,
$j\in\w$. Each $\W_n$ is a point finite family since both $\B_n$ and
$\C_n$ are point finite and so both $\{B_j\}_{j\in\w}$ and
$\{C_j\}_{j\in\w}$ are infinite. 
Since $\B$ and $\C$ are sharp bases, this implies that 
$\{\bigcap_{j\le n}B_j\times C_j:n\in\w\}$ is a local base at $(x,r)$ and
$\W$ is a sharp base as required.
\end{proof}

Ponomarev, see \cite{G}, characterized those spaces with a point countable
base as precisely the open $s$-images of metric spaces (a map is an
$s$-map if it has separable fibres). There is a similar characterization
for sharp bases.

\begin{thm} 
A space $X$ has a sharp base if and only if there is a metric space $M$
with a base $\B$ and a continuous open mapping $f:M\to X$ such that,
whenever $x\in X$ and $\{B_n\in\B:n\in\w\}$ is a pairwise distinct
collection, if $f^{-1}(x)\cap B_n\neq\nowt$ for each $n\in\w$, then there
exists $n_0$ such that for each $y\in X$, if $f^{-1}(y)\cap B_j\neq\nowt$,
for each $j\le n_0$, then $f^{-1}(y)\cap B_0\neq\nowt$.
\end{thm}

\begin{proof} 
Suppose that $\G$ is a sharp base for the space $X$. Let
$$
M=\{(G_n)\in\G^\w:x\in\bigcap_{n\in\w}G_n\text{ for some $x\in X$}\}
$$
be the subspace of the Baire metric space $\G^\w$, with metric 
$d((G_n),(H_n))=1/2^k$ where $k$ is least such that $G_n\neq H_n$.
Let $f:M\to X$ be defined letting $f((G_n))$ be the unique element of
$\bigcap_{n\in\w}G_n$ and let $\B$ be the base for $M$ consisting of all
$1/2^n$-balls about points of $M$.
Then $f$ is easily seen to be a continuous, open mapping onto $X$ and the 
condition on $\B$ in the statement of the theorem is merely a translation
of the fact that $\G$ is a sharp base.
\end{proof}

It is clear from the proof that, in the statement of the theorem, we can
take $\B$ to be the collection of $1/2^n$ balls for any $n$ rather than a
base for $M$. Since a space with a sharp base has a point countable sharp
base, we can also assume that the map in the statement of the theorem is
an $s$-map. However, it is not immediately clear that we can prove that a
space with a sharp base has a point countable base directly from the
theorem.

We conclude with some open problems. Since every collectionwise normal
Moore space is metrizable, the following is a natural and intriguing
question.

\begin{qu} 
Is every collectionwise normal space with a sharp base metrizable?
\end{qu}

Example 4 of \cite{AA} shows that weakly developable, collectionwise 
normal spaces do not have to be metrizable and the Heath V-space over a
Q-set is an example of a normal space with a uniform base that is not
metrizable. On the other hand, the answer is \lq yes\rq\ if the space is
also submetacompact (since it is then a Moore space) or a strict p-space.
We might also ask whether a perfect, collectionwise normal space
with a sharp base is metrizable.
It is interesting to note that it is not known whether a collectionwise
normal space with a point countable base need be paracompact.

Since the Heath V-space over a $\Delta$-set is countably paracompact but
not normal \cite{kn}, at least consistently a countably paracompact,
(Moore) space with a sharp base need not be normal. What about the
converse?

\begin{qu} 
Is there a Dowker space with a sharp base?
\end{qu}

\begin{qu} 
Is every perfectly regular space with a sharp base developable? 
Is every normal space with a sharp base developable? 
Is every perfectly regular, pseudocompact space with a sharp base
metrizable?
\end{qu}

Not every Moore space with a weakly uniform base has a uniform base (see 
\cite{AA}) so we ask:

\begin{qu} 
Does every Moore space with a sharp base have a uniform base?
\end{qu}

Every pseudocompact space with a $G_\dlt$-diagonal is \v Cech complete
\cite{A}, and every pseudocompact Moore space with a sharp base is
metrizable.

\begin{qu} Is every \v Cech complete Moore space with a sharp base
metrizable? What about Baire instead of \v Cech complete?
\end{qu}

\begin{qu} If $X\times[0,1]$ has a sharp base, does $X$ have a
$\sgm$-point finite sharp base?
\end{qu}

\begin{qu}
Does the image (or pre-image) of a space with a sharp base under a
perfect map (closed and open map, open map with compact, countable or
finite fibres) have a sharp base?
\end{qu}

\subsection*{Acknowledgements} 
Much of this research took place whilst the first author held a
Universitas 21 Traveling Fellowship at the Department of Mathematics at
the University of Auckland. 
He would like to thank the department for their hospitality and both the
department and U21 for their financial support.


\begin{thebibliography}{10}

\bibitem{al}
P.~Aleksandrov, \emph{On the metrisation of topological spaces}, Bull. Acad.
  Polon. Sci. S\'er. Sci. Math. Astr. Phys. \textbf{8} (1960), 135--140. \MR{22
  \#5024}

\bibitem{AA}
B.~Alleche, A.~V. Arhangel$'$ski{\u\i}, and J.~Calbrix, \emph{Weak developments
  and metrization}, Topology Appl. \textbf{100} (2000), no.~1, 23--38, Special
  issue in honor of Howard H. Wicke. \MR{2001g:54029}

\bibitem{arh}
A.~V. Arhangel$'$ski{\u\i}, \emph{Some metrization theorems}, Uspehi Mat. Nauk
  \textbf{18} (1963), no.~5 (113), 139--145. \MR{27 \#6242}

\bibitem{A}
A.~V. Arhangel$'$ski{\u\i}, W.~Just, E.~A. Rezniczenko, and P.~J. Szeptycki,
  \emph{Sharp bases and weakly uniform bases versus point-countable bases},
  Topology Appl. \textbf{100} (2000), no.~1, 39--46, Special issue in honor of
  Howard H. Wicke. \MR{2001f:54027}

\bibitem{br}
Z.~Balogh and M.~E. Rudin, \emph{Monotone normality}, Topology Appl.
  \textbf{47} (1992), no.~2, 115--127. \MR{94b:54065}

\bibitem{BL}
Harold Bennett and David Lutzer, \emph{Ordered spaces with special bases},
  Fund. Math. \textbf{158} (1998), no.~3, 289--299. \MR{99i:54044}

\bibitem{BL2}
\bysame, \emph{Continuous separating families in ordered spaces and strong base
  conditions}, Topology Appl. \textbf{119} (2002), 305--314.

\bibitem{gi}
Raymond~F. Gittings, \emph{Concerning quasi-complete spaces}, General Topology
  and Appl. \textbf{6} (1976), no.~1, 73--89. \MR{52 \#11855}

\bibitem{gkm}
Chris Good, Robin~W. Knight, and Abdul~M. Mohamad, \emph{On the metrizability
  of spaces with a sharp base}, To appear in Topology Appl. PII:
  S0166-8641(01)00300-5, 2002.

\bibitem{G}
Gary Gruenhage, \emph{Generalized metric spaces}, Handbook of set-theoretic
  topology, North-Holland, Amsterdam, 1984, pp.~423--501. \MR{86h:54038}

\bibitem{HL}
R.~W. Heath and W.~F. Lindgren, \emph{Weakly uniform bases}, Houston J. Math.
  \textbf{2} (1976), no.~1, 85--90. \MR{52 \#15364}

\bibitem{H}
R.~E. Hodel, \emph{Metrizability of topological spaces}, Pacific J. Math.
  \textbf{55} (1974), 441--459. \MR{51 \#6747}

\bibitem{kn}
Robin~W. Knight, \emph{$\delta$-sets}, Trans. Amer. Math. Soc. \textbf{339}
  (1993), no.~1, 45--60. \MR{94a:54016}

\bibitem{Mc}
William~G. McArthur, \emph{${G}\sb{\delta }$-diagonals and metrization
  theorems}, Pacific J. Math. \textbf{44} (1973), 613--617. \MR{47 \#5835}

\bibitem{per}
S.~A. Peregudov, \emph{On pseudocompactness and other covering properties},
  Questions Answers Gen. Topology \textbf{17} (1999), no.~2, 153--155.
  \MR{2000h:54033}

\bibitem{sc}
Brian~M. Scott, \emph{Pseudocompact, metacompact spaces are compact}, Topology
  Proc. \textbf{4} (1979), no.~2, 577--587 (1980), The Proceedings of the 1979
  Topology Conference (Ohio Univ., Athens, Ohio, 1979). \MR{81m:54034}

\bibitem{sha}
D.~B. Shakhmatov, \emph{Pseudocompact spaces with a point-countable base},
  Dokl. Akad. Nauk SSSR \textbf{279} (1984), no.~4, 825--829, English
  translation: Soviet Math. Dokl. 30 (1984), no. 3, 747--751. \MR{86a:54029}

\bibitem{us}
V.~V. Uspenski{\u\i}, \emph{Pseudocompact spaces with a $\sigma$-point-finite
  base are metrizable}, Comment. Math. Univ. Carolin. \textbf{25} (1984),
  no.~2, 261--264. \MR{87f:54035}

\bibitem{w}
W.~Stephen Watson, \emph{Pseudocompact metacompact spaces are compact}, Proc.
  Amer. Math. Soc. \textbf{81} (1981), no.~1, 151--152. \MR{81j:54042}

\bibitem{wat}
\bysame, \emph{A pseudocompact meta-{L}indel\"of space which is not compact},
  Topology Appl. \textbf{20} (1985), no.~3, 237--243. \MR{87d:54042}

\end{thebibliography}
\providecommand{\bysame}{\leavevmode\hbox to3em{\hrulefill}\thinspace}
\providecommand{\MR}{\relax\ifhmode\unskip\space\fi MR }
\providecommand{\MRhref}[2]{%
  \href{http://www.ams.org/mathscinet-getitem?mr=#1}{#2}
}
\providecommand{\href}[2]{#2}

\end{document}